\documentstyle[12pt,epic,eepic]{ioplppt}
\input epsf.tex
\begin{document}

\jl{1}

\title{Ribbon tableaux and $q$-analogues of fusion rules\\
in WZW conformal field theories}[]

\author{Omar Foda\dag,  Bernard Leclerc\ddag\,, Masato Okado\S\
and Jean-Yves Thibon$\Vert$}

\address{\dag Department of Mathematics, The University of Melbourne,
Parkville, Victoria 3052, Australia}

\address{\ddag D\'epartement de Math\'ematiques,
Universit\'e de Caen, Campus II,
Boulevard Mar\'echal Juin,
BP 5186, 14032 Caen cedex, France}

\address{\S Department of Mathematical Science, Faculty of Engineering Science,
Osaka University, Osaka 560-8531, Japan}

\address{$\Vert$ Institut Gaspard Monge, Universit\'e de Marne-la-Vall\'ee,
Cit\'e Descartes, 5 Boulevard Descartes, Champs-sur-Marne,
77454 Marne-la-Vall\'ee cedex, France}


\begin{abstract}
Starting from known
$q$-analogues of ordinary $SU(n)$ tensor products multiplicities,
we  introduce $q$-analogues of the  fusion coefficients of the
WZW conformal field theories associated with $SU(n)$. We conjecture
combinatorial interpretations of these polynomials,
which can be proved in special cases. This allows us to derive in a simple way
various kinds of  branching functions, the simplest ones being
the characters of the minimal unitary series of the Virasoro algebra.
We also obtain $q$-analogues of the dimensions of spaces
of nonabelian theta functions.
\end{abstract}

\font\twelvegoth=eufm10 at 12pt
\font\tengoth=eufm10
\font\sevengoth=eufm7
\font\fivegoth=eufm5

\newfam\gothfam
\textfont\gothfam=\twelvegoth
\scriptfont\gothfam=\sevengoth
\scriptscriptfont\gothfam=\fivegoth
\def\goth{\fam\gothfam\tengoth}

\def\SG{{\goth S}}
\def\SP{{\rm \widehat{S}}}
\def\sym{Sym}
\def\A{{\cal A}}
\def\L{{\cal L}}
\def\P{{\cal P}}
\def\bs{\bar{s}}
\def\bk{\bar{k}}
\def\C{{\bf C}}
\def\Z{{\bf Z}}
\def\Q{{\bf Q}}
\def\N{{\bf N}}
\def\W{\hat{W}}
\def\Hom{{\rm Hom}}
\def\STab{{\rm STab}}
\def\slchap{\widehat{sl}}
\def\Des{{\rm Des\,}}
\def\maj{{\rm maj\,}}
\def\mod{{\ \rm mod\ }}
\def\<{\langle}
\def\>{\rangle}
\def\pr#1{\langle #1 \rangle}
\def\prh#1{\pr{\widehat{#1}}}
\def\vac{|0\>}
\def\DP{{\rm DP}}
\def\DPR{{\rm DPR}}
\def\SR{{\rm SR}}
\def\ket#1{|#1\>}
\def\Coinv{{\cal H}}
\def\ch{{\rm ch\,}}
\def\Tab{{\rm Tab\,}}
\def\e{{\bf e}}
\def\lev{{\rm lev\,}}
\def\wt{{\rm wt\,}}
\def\ft{\tilde{f}}
\def\et{\tilde{e}}
\def\ie{{\it i.e.\,}}
\def\ot{\otimes}

\def\shuff#1#2{\mathbin{
      \hbox{\vbox{
        \hbox{\vrule
              \hskip#2
              \vrule height#1 width 0pt
               }%
        \hrule}%
             \vbox{
        \hbox{\vrule
              \hskip#2
              \vrule height#1 width 0pt
               \vrule }%
        \hrule}%
}}}

\def\shuffle{{\mathchoice{\shuff{7pt}{3.5pt}}%
                        {\shuff{6pt}{3pt}}%
                        {\shuff{4pt}{2pt}}%
                        {\shuff{3pt}{1.5pt}}}}%
\def\shuf{\,\shuffle\,}
\def\qshuf{\shuf_q}

\def\qbin#1#2{\left[\matrix{#1 \cr #2}\right]_q}

\newtheorem{example}{Example}[section]
\newtheorem{theorem}[example]{Theorem}
\newtheorem{corollary}[example]{Corollary}
\newtheorem{definition}[example]{Definition}
\newtheorem{proposition}[example]{Proposition}
\newtheorem{algorithm}[example]{Algorithm}
\newtheorem{lemma}[example]{Lemma}
\newtheorem{conjecture}[example]{Conjecture}

\font\FontSeteight=msbm8
\font\FontSetnine=msbm9
\font\FontSetten=msbm10
\font\FontSettwelve=msbm10 scaled 1200
\newdimen\Squaresize \Squaresize=15pt
\newdimen\Thickness \Thickness=0.5pt

\def\Square#1{\hbox{\vrule width \Thickness
   \vbox to \Squaresize{\hrule height \Thickness\vss
      \hbox to \Squaresize{\hss#1\hss}
   \vss\hrule height\Thickness}
\unskip\vrule width \Thickness}
\kern-\Thickness}

\def\Vsquare#1{\vbox{\Square{$#1$}}\kern-\Thickness}
\def\blk{\omit\hskip\Squaresize}

\def\young#1{
\vbox{\smallskip\offinterlineskip
\halign{&\Vsquare{##}\cr #1}}}

\def\thisbox#1{\kern-.09ex\fbox{#1}}
\def\downbox#1{\lower1.200em\hbox{#1}}

\section{Introduction}

The Littlewood-Richardson coefficient
$c^\lambda_{\mu^{(1)},\ldots,\mu^{(r)}}$,
where $\lambda,\mu^{(1)},\ldots,\mu^{(r)}$ are partitions
of $N,m_1,\ldots,m_r$ ($N=m_1+\cdots+m_r$) is defined as the multiplicity
of the irreducible representation $V_\lambda$ of $U(n)$
(or of $gl(n,\C)$) in the tensor product
$V_{\mu^{(1)}}\otimes\cdots\otimes V_{\mu^{(r)}}$,
or equivalently, as the coefficient of the Schur function $s_\lambda$
in the product $s_{\mu^{(1)}}\cdots s_{\mu^{(r)}}$. Also, by
Schur-Weyl duality, $c^\lambda_{\mu^{(1)},\ldots,\mu^{(r)}}$ is
the multiplicity of the irreducible representation $S_\lambda$
of $\SG_N$ in the representation induced from the representation
$S_{\mu^{(1)}}\otimes\cdots\otimes S_{\mu^{(r)}}$ of
the Young subgroup $\SG_{m_1}\times\cdots\times \SG_{m_r}$.

In recent years, two kinds of generalizations of these numbers have
been considered:

\medskip
1) Polynomials $c^\lambda_{\mu^{(1)},\ldots,\mu^{(r)}}(q)$ with
integer coefficients, reducing to
Littlewood-Richardson multiplicites for $q=1$. As a special case,
one finds the Kostka-Foulkes polynomials, which can also be
interpreted as $q$-analogues of weight multiplicities. The most
general $q$-analogues are defined in terms of certain generalized
Young tableaux, called {\it ribbon tableaux}  \cite{LLTrib}.
These polynomials are known to be related to the quantum affine algebras
$U_q(\slchap_n)$, and have been recently identified as a special
family of Kazhdan-Lusztig polynomials for affine symmetric groups
\cite{LT98}. Their specialization at roots of unity are believed
to be relevant to the calculation of certain plethysms \cite{LLTrib}.
Another definition of $q$-analogues of LR coefficients
appears in \cite{SW}. It seems to coincide with ours in most cases
but the reason for this is still unclear.

\medskip
2) Restricted Littlewood-Richardson coefficients
$\bar c^\lambda_{\mu^{(1)},\ldots,\mu^{(r)}}$, which are the structure
constants of the fusion algebras of Wess-Zumino-Witten conformal
field theories associated to $(\slchap_n)_l$.
These numbers can be computed as alternating sums of ordinary
Littlewood-Richardson coefficients (the Kac-Walton formula \cite{Wal1,Kac,Cum}).
They are known to be nonnegative, but a combinatorial interpretation
is available only in special cases. 

\bigskip
The aim of this paper is to propose a common generalization
of 1) and 2), that is, to define $q$-analogues of the restricted
Littlewood-Richardson coefficients in terms of ribbon tableaux.
Our definition proceeds by introducing  appropriate powers of $q$
in the Kac-Walton formula, so that the replacement of the ordinary
multiplicities by their $q$-analogues yields a nonnegative polynomial.
We are able to prove this only in special cases, but the simplest
examples are already of interest, since they lead to simple formulae
for various kinds of branching functions, such as the characters
of the representations of unitary minimal series of the Virasoro algebra.
We also obtain $q$-analogues of the dimensions of the spaces of
nonabelian theta functions.

\section{Fusion rules and vertex operators}

Let $\Lambda_0,\Lambda_1,\ldots,\Lambda_{n-1}$ be the fundamental
weights of the affine Lie algebra $\slchap_n$ of type $A^{(1)}_{n-1}$,
and let $\delta$ be the null root. The Chevalley generators are
denoted by $e_i,h_i,f_i$, $i=0,\ldots,n-1$ and the degree generator
by $d$. The weight lattice is 
$P=\Z\Lambda_0\oplus\cdots\oplus\Z\Lambda_{n-1}\oplus\Z\delta$.
The set of dominant integral weights is
$P^+=\N\Lambda_0\oplus\cdots\oplus\N\Lambda_{n-1}\oplus\Z\delta$.
The Weyl vector is $\rho=\Lambda_0+\cdots +\Lambda_{n-1}$.
The level of a weight $\Lambda=\sum_{i=0}^{n-1}l_i\Lambda_i+z\delta\in P^+$ 
is $\lev (\Lambda)=\sum_{i=0}^{n-1}l_i$.
We denote by $P_l^+$ the subset of level $l$ dominant integral
weights of the form $\Lambda=\sum_{i=0}^{n-1}l_i\Lambda_i$
(no $\delta$)
with $l_i\in\N$ and $l_0+\cdots+l_{n-1}=l$.
For each $\lambda\in P^+$ there is a unique integrable highest
weight module $L(\Lambda)$ with highest weight $\Lambda$.
Any weight $\mu\in P^+$ of the form
$\mu=m_1\Lambda_1+\cdots+m_{n-1}\Lambda_{n-1}$ can be interpreted
as a dominant integral weight of the finite dimensional Lie algebra
$sl_n\subset\slchap_n$. For any choice of a complex number
$z$, the irreducible representation $V_\mu$ of $sl_n$
can be extended to a level 0 representation $V_\mu(z)$ of
${\slchap_n}'$ (evaluation modules).

The classical part $\lambda=\bar\Lambda$ of a weight $\Lambda\in P$
is defined by means of the linear map $\bar\Lambda_i=\Lambda_i-\Lambda_0$,
$\bar\delta=0$.
If we fix a level $l\ge 1$, a dominant integral weight
$\Lambda\in P_l^+$ is  specified by its classical part $\lambda=\bar\Lambda$,
which can be identified with a partition whose Young diagram is contained
in the rectangle $(l^{n-1})$, the fundamental weight $\Lambda_i$ being
represented by a column of height $i$. Conversely, we set $\wt (\lambda)=\Lambda$.

For such partitions $\lambda,\mu^{(1)},\ldots,\mu^{(r)}$, let
$\bar c^\lambda_{\mu^{(1)},\ldots,\mu^{(r)}}$ denote the dimension
of the space of intertwining operators
\begin{equation} \label{eq:dim-of-VO}
\Hom_{{\slchap_n}'}(L(\Lambda),\, L(l\Lambda_0)
\otimes V_{\mu^{(1)}}(z_1)\otimes\cdots\otimes V_{\mu^{(r)}}(z_r) )
\end{equation}
(called chiral vertex operators) where $(z_1\ldots,z_r)\in\C^r$
is generic.
It is known that the numbers $\bar c^\lambda_{\mu^{(1)},\ldots,\mu^{(r)}}$
are the structure constants of an associative algebra
${\cal F}^{(n,l)}$, the fusion algebra of the Wess-Zumino-Witten model
associated to $\slchap_n$ at level $l$
(see \cite{TUY}). This algebra has several
(non obviously equivalent) interpretations. The most elementary one,
on which we will rely in this paper, is due to Goodman and Wenzl \cite{GW}.

Let $H_N(q)$ be the Hecke algebra of type $A_{N-1}$, \ie the
$\C[q,q^{-1}]$-algebra generated by elements $T_1,\ldots,T_{N-1}$
verifying the braid relations together with $T_i^2=(q-1)T_i+q$.
This algebra is a $q$-analogue of the group algebra of the symmetric group
$\SG_N$, and there exist for it $q$-analogues of the various realizations
of the irreducible representations of the symmetric group. In particular,
there is a $q$-analogue $\pi_\lambda^q$ of Young's orthogonal form.
The representation space $V_\lambda(q)$ has an orthonormal basis
$|t\>$ labelled by standard Young tableaux of shape $\lambda$.

Recall that a standard tableau $t$ can be interpreted as a chain
of partitions $\lambda^{(1)}\subset \lambda^{(2)}\subset\ldots\subset\lambda^{(N)}
=\lambda$, where $\lambda^{(k)}$ is obtained by adding to the diagram
of $\lambda^{(k-1)}$ the box containing the entry $k$ in $t$.
Let us say that a partition $\lambda$ is $(n,l)$-restricted if
it has at most $n$ parts and $\lambda_1-\lambda_n\le l$. Let
$\Pi^{(n,l)}$ be the set of such partitions. The tableau $t$
is said to be $(n,l)$-restricted if all the intermediate partitions
$\lambda^{(k)}\in\Pi^{(n,l)}$. Denote by $\STab^{(n,l)}(\lambda)$
the set of $(n,l)$-restricted tableaux of shape $\lambda$.

Suppose now that $q=\zeta$, a primitive $L$th root of unity, where
$L=n+l$. Then, $H_N(\zeta)$ is not semisimple for $N\ge L$, and
$V_\lambda(\zeta)$ is not irreducible, nor even semisimple in general.
Moreover only the integral form of the irreducible representations
can be specialized at $q=\zeta$.
For $\lambda\in\Pi^{(n,l)}$, Wenzl \cite{We} showed how to construct
{}from the orthogonal form
$V_\lambda(\zeta)$ an irreducible representation $D_\lambda$.
The point is that the matrix elements $\<t'|\pi_\lambda^q(T_w)|t''\>$
have no pole at $q=\zeta$ as soon as $t'$ and $t''$ are $(n,l)$-retricted.
Moreover, the representation matrices restricted to the subspace
\begin{equation}
V^{(n,l)}_\lambda(\zeta) = \bigoplus_{t\in\STab^{(n,l)}(\lambda)}\C |t\>
\end{equation}
define an irreducible representation $D_\lambda$ of $H_N(\zeta)$.
Denote by $\pi_\lambda^{(n,l)}$ the representation on this space, and
let
\begin{equation}
H_N^{(n,l)} = \bigoplus_{\lambda\in \Pi_N^{(n,l)}} 
\pi^{(n,l)}_\lambda(H_N(\zeta))\,.
\end{equation}
Then, $H_N^{(n,l)}$ is a semi-simple quotient of $H_N(\zeta)$,
and its irreducible representations are exactly the $\pi^{(n,l)}_\lambda$.

Let ${\cal R}_N^{(n,l)}=R(H_N^{(n,l)})$ be the Grothendieck group
of $H_N^{(n,l)}$, \ie the free abelian group
generated by the isomorphism classes $[D_\lambda]$ of irreducible
representations, with addition corresponding to direct sum. The sum
\begin{equation}
{\cal R}^{(n,l)}=\bigoplus_{N\ge 0}{\cal R}_N^{(n,l)}
\end{equation}
can be endowed with a ring structure by setting
$[D_\lambda]\cdot [D_\mu] =
\left[ D_\lambda\hat\otimes D_\mu\right] $
where $\hat\otimes$ is the external tensor product, obtained by inducing
the $H_N\otimes H_M$-module $D_\lambda\otimes D_\mu$ to $H_{N+M}$.
Then \cite{GW}, ${\cal R}^{(n,l)}$ is isomorphic to a quotient
$Sym(n)/{\cal J}^{n,l}$ of the ring of symmetric polynomials in
$n$ variables $Sym(n)=\Z[x_1,\ldots,x_n]^{\SG_n}$, the ideal
${\cal J}^{n,l}$ being generated by the Schur functions
$s_\lambda$ labelled by partitions such that $\lambda_1-\lambda_n=l+1$.
The fusion algebra ${\cal F}^{(n,l)}$ is isomorphic to the quotient
of ${\cal R}^{(n,l)}$ by the single relation $s_{(1^n)}\equiv 1$.

\section{Crystal base and $q$-vertex operators}

If we consider $q$ as an indeterminate, we can replace $\slchap_n'$
in (\ref{eq:dim-of-VO}) by $U_q'(\slchap_n)$. According to \cite{FR,DJO}
the space of $q$-vertex operators $\Hom_{U_q'(\slchap_n)}(L(\Lambda),
L(l\Lambda_0)\otimes V_{\mu^{(1)}}(z_1)\otimes\cdots\otimes V_{\mu^{(r)}}(z_r))$
is identified with the following vector space:
\begin{equation} \label{eq:sp-of-VO}
\bigoplus\bigotimes_{j=1}^r\Q(q)\< v\in V_{\mu^{(j)}}\mid
\wt v=\Lambda^{(j)}-\Lambda^{(j-1)},e_i^{\<h_i,\Lambda^{(j-1)}\>+1}
v=0\,\forall i\>,
\end{equation}
where the direct sum is taken over all sequences 
$(\Lambda^{(0)},\ldots,\Lambda^{(r)})\in(P^+_l)^{r+1}$ such that 
$\Lambda^{(0)}=l\Lambda_0,\Lambda^{(r)}=\Lambda$.

We would like to relate our restricted LR coefficient to the crystal base
theory introduced by Kashiwara \cite{Ka}. For this purpose we prepare some 
notations. Any integrable $U_q(\slchap_n)$-module $L(\Lambda)$ has a crystal
base $(\L(\Lambda),B(\Lambda))$ \cite{Ka}. But the finite-dimensional 
$U_q'(\slchap_n)$-module $V_\lambda$ does not necessarily have a crystal base.
{}From this reason we restrict ourselves later in this section
to the cases when all $\mu^{(j)}$'s are of rectangular shape, and let 
$(\L_{\mu^{(j)}},B_{\mu^{(j)}})$
be the crystal base of $V_{\mu^{(j)}}$. As is well known, modified Chevalley 
generators (or Kashiwara operators) $\et_i,\ft_i$ act on $B(\Lambda)$ or 
$B_{\mu^{(j)}}$.

Let $\Phi(z)$ be an appropriately normalized $q$-vertex operator from
$L(\Lambda')$ to $L(\Lambda)\otimes V_\lambda(z)$.
It is known \cite{DJO} that it preserves the crystal lattice, \ie
\[
\Phi(z)\L(\Lambda')\subset\L(\Lambda)\otimes\L_\lambda(z).
\]
Therefore, counting the dimensionality of the space of vertex operators 
is reduced to the combinatorics of crystals. Let us define the set of restricted
paths.
\[
\P^\lambda_{\mu^{(1)},\ldots,\mu^{(r)}}
=\left\{p=b_1\otimes\cdots\otimes b_r\left|
\begin{array}{l}
b_j\in B_{\mu^{(j)}}, \wt b_j=\Lambda^{(j)}-\Lambda^{(j-1)},\\
\et_i^{\<h_i,\Lambda^{(j-1)}\>+1}b_j=0\,\forall j=1,\ldots,r
\end{array}
\right.\right\},
\]
where $\Lambda^{(j)}$'s are as in (\ref{eq:sp-of-VO}). Thus we have 
\[
\bar c^\lambda_{\mu^{(1)},\ldots,\mu^{(r)}}
=\left|\P^\lambda_{\mu^{(1)},\ldots,\mu^{(r)}}\right|.
\]

We shall introduce a statistic on the set of restricted paths.
First we introduce a $\Z$-valued function $H_{\mu^{(j)}\mu^{(k)}}$
on $B_{\mu^{(j)}}\otimes B_{\mu^{(k)}}$ called the energy function.
It is determined from the combinatorics of crystal graphs.
For an element $p=b_1\ot\cdots\ot b_r$ of 
$B_{\mu^{(1)}}\ot\cdots\ot B_{\mu^{(r)}}$, we define the energy of $p$ by
\[
E(p)=\sum_{1\le j<k\le r}H_{\mu^{(j)}\mu^{(k)}}(b_j\ot b^{(j+1)}_k),
\]
where $b^{(j+1)}_k$ is determined under some isomorphism of crystals
(see \cite{NY,SWa,Sh} for details).
Using this statistic we define a $q$-analogue 
$\bar c^\lambda_{\mu^{(1)},\ldots,\mu^{(r)}}(q)$ of 
$\bar c^\lambda_{\mu^{(1)},\ldots,\mu^{(r)}}$ by
\begin{equation}\label{LRpath}
\bar c^\lambda_{\mu^{(1)},\ldots,\mu^{(r)}}(q)
=\sum_{p\in{\P^\lambda_{\mu^{(1)},\ldots,\mu^{(r)}}}}q^{E(p)}.
\end{equation}

\section{Ribbon tableaux and $q$-analogues of LR coefficients}

Recall that a Schur function $s_\lambda(X)$ can be expressed as
a sum over semi-standard Young tableaux $t$ of shape $\lambda$
\begin{equation}
s_\lambda(X)=\sum_{t\in\Tab(\lambda)}X^t 
\end{equation}
where $X^t=\prod_i x_i^{m_i}$, $m_i$ being the number of occurences
of the integer $i$ in $t$.
Therefore, a product of $r$ Schur functions $s_{\mu^{(i)}}$
is a sum over $r$-tuples of tableaux
\begin{equation}
s_{\mu^{(1)}}s_{\mu^{(r)}}\cdots s_{\mu^{(r)}}
=
\sum_{(t_1,\ldots,t_r)}X^{t_1}X^{t_2}\cdots X^{t_r}\,.
\end{equation}
The point of this trivial remark is that $r$-tuples of tableaux
are in one-to-one correspondence with a certain kind of generalized
Young tableaux, the so-called $r$-ribbon tableaux, on which some
extra combinatorial information can be read.

The notion of ribbon tableau is derived from the $r$-core and
$r$-quotient algorithms, originating from the modular representation theory
of symmetric groups ({\it cf.} \cite{JK}). An $r$-ribbon $R$ is a connected skew
Young diagram of $r$ boxes, not containing a $2\times 2$ square.
Its height $h(R)$ is the number of rows occupied by the diagram.
For example,
%
\begin{center}
\leavevmode
\epsfxsize =4cm
\epsffile{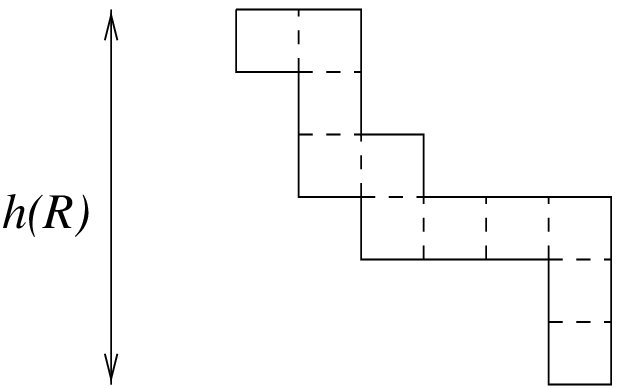}
\end{center}
%
is an 11-ribbon of height $h(R)=6$. A partition is said to be an
$r$-core if it is not possible to peel off an $r$-ribbon from
its diagram. The $r$-core $\lambda_{(r)}$ 
of a partition $\lambda$ is the partition
obtained by removing from its diagram a maximal number of $r$-ribbons.
This can usually be done in several ways, but  all
of them lead to the same result. The $r$-core play the role of the remainder
in a kind of Euclidean division of partitions. The role of the quotient
is played by an $r$-tuple of partitions,
the $r$-quotient  $(\lambda^{(0)},\ldots,\lambda^{(r-1)})$,
satisfying
\begin{equation}
|\lambda|=|\lambda_{(r)}|+r\sum_{i=0}^{r-1}|\lambda^{(i)}|\,.
\end{equation}
Details  can be found in \cite{JK}. This algorithm
is revertible, and it provides a one-to-one correspondence between
partitions with fixed $r$-core $\nu$ and $r$-tuples of partitions.
In particular, any $r$-tuple of partitions can  be interpreted
as a single partition with empty $r$-core. 

In the same way as standard tableaux can be regarded as chains of
partitions whose consecutive terms differ by exactly one box, semi-standard
tableaux can be interpreted as chains whose consecutive terms differ
by horizontal strips of boxes. Applying the inverse $r$-quotient algorithm
to an $r$-tuple of semi-standard tableaux, interpreted as a chain of $r$-tuples
of partitions, one obtains a chain of partitions without $r$-core,
in which  the skew diagrams formed from two consecutive terms $k-1,k$ have
the property of being tilable by $r$-ribbons in exactly one way.
These ribbons can be labelled by $k$, and the chain can be represented
by an $r$-ribbon tableau, which is a tiling of a Young diagram by
labelled $r$-ribbons, satisfying some order conditions (precisely,
if we define the root of a ribbon as its rightmost lowest cell,
the root of a ribbon labelled $k$ should not lie above any ribbon
labelled $j$ for $j\ge k$). The weight of a ribbon tableau is defined
as for ordinary tableaux, so that the correspondence between
$r$-tuples of ordinary tableaux and ribbon tableaux is weight
preserving. The spin $s(R)$ of a ribbon $R$ is ${1\over 2}(h(R)-1)$, and
the spin $s(T)$ of a ribbon tableau $T$ is the sum of the spins of its ribbons.
For example,
\begin{center}
\leavevmode
\epsfxsize =4cm
\epsffile{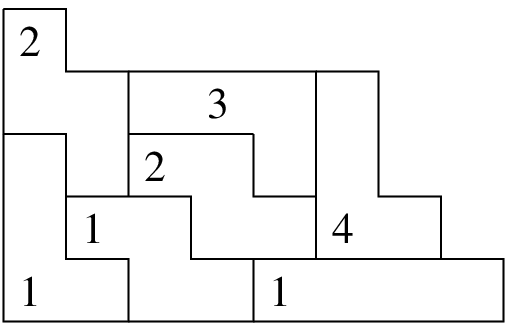}
\end{center}
%
is a 4-ribbon tableau of shape $(8,7,6,6,1)$, weight $(3,2,1,1)$ and spin $9/2$.

Let $T$ be a ribbon tableau of shape $\lambda$. Let $s_r^*(\lambda)$
be the maximal spin of an $r$-ribbon tableau of shape $\lambda$, and
define the cospin $\tilde s(T)$ as the difference $s_r^*(\lambda)-s(T)$.
It can be shown that this is always an integer. For a
partition $\lambda$ without $r$-core, let $\Tab_r(\lambda)$ be the
set of semi-standard $r$-ribbon tableaux of shape $\lambda$. Then,
the polynomial
\begin{equation}
\tilde G_\lambda^{(r)}(X;q)=
\sum_{T\in\Tab_r(\lambda)} q^{s(T)} X^T
\end{equation}
is a $q$-analogue of the product of Schur functions
$s_{\lambda^{(0)}}s_{\lambda^{(1)}}\cdots s_{\lambda^{(r-1)}}$.
It can be shown that it is actually a symmetric polynomial, though
this is not obvious from the definition. The polynomials
\begin{equation}
c^\nu_{\lambda^{(0)},\ldots,\lambda^{(r-1)}}(q)
= \<s_\nu,\, \tilde G_\lambda^{(r)}(X;q)\>
\end{equation}
are therefore $q$-analogues of Littlewood-Richardson coefficients.
It is conjectured that their coefficients are non-negative (this is
also unclear from their definition). For $r=2$ this has been proved
by a combinatorial argument \cite{CL}. This is also known when
all the $\lambda^{(i)}$ are row (resp. column) diagrams. 
For the general
case, the positivity should follow from a recent expression of
these polynomials as parabolic Kazhdan-Lusztig polynomials of the
affine symmetric group \cite{LT98}, but the relevant geometrical
interpretation does not seem to be available in the literature.

An interesting particular case is obtained when $\lambda$ is of
the form $\lambda=r\mu=(r\mu_1,\ldots,r\mu_m)$, $\mu$ being an arbitrary
partition. Then $\lambda$ has no $r$-core, and one can define
\begin{eqnarray}
\tilde H_\mu^{(r)}(X;q) &= \tilde G_{r\mu}^{(r)}(X;q)\,, \\
H_\mu^{(r)}(X;q) &= q^{s_r^*(r\mu)}\tilde H_\mu^{(r)} (X;q^{-1})\,.
\end{eqnarray}
For $r$ sufficiently large ($r\ge m=\ell(\mu)$), it is known that
$H_\mu^{(r)}$ is equal to the Hall-Littlewood function
\begin{equation}
Q'_\mu(X;q)=\sum_\lambda K_{\lambda\mu}(q) s_\lambda
\end{equation}
where the $K_{\lambda\mu}(q)$ are the Kostka-Foulkes polynomials.
Moreover, it is conjectured that the difference between two
consecutive $H$-functions $H_\mu^{(j+1)}-H_\mu^{(j)}$ is nonnegative
on the Schur basis.

\section{The $(n,l)$-restricted $q$-analogues}

One can compute in fusion algebras ${\cal F}^{(n,l)}$
by an algorithm due (independently) to Kac and Walton. In the case of
$\slchap_n$, this algorithm turns out to be identical
to the one devised by Goodmann and Wenzl for the algebras
${\cal R}^{(n,l)}$, and it is this coincidence which proves
the relation between these algebras.

The Goodman-Wenzl algorithm for computing in ${\cal R}^{(n,l)}$ can
be described as follows. To compute the product
$\bar{f}_1\cdots \bar{f}_r$ of elements $\bar{f}_i\in {\cal R}^{(n,l)}$,
one first evaluates the corresponding product of symmetric functions
$$
f_1\cdots f_r = \sum_\lambda c_\lambda s_\lambda
$$
in terms of Schur functions. Next, one replaces each $s_\lambda$
by its class $\bs_\lambda$, taking into account the equivalences
\begin{equation}
\bs_\lambda =
\cases{\varepsilon(w)\bar s_{w\circ\lambda} & if there is a $w\in\W$ such that 
                                          $w\circ\lambda\in\Pi^{(n,l)}$\\
                         0 & otherwise\\}
\end{equation}
where the action of the affine Weyl group $\W=\<\sigma_0,\ldots,\sigma_{n-1}\>$ 
of type $A^{(1)}_{n-1}$
on $\Z^n$ is
defined by
\begin{equation}
\sigma_i(\lambda) = \cases{(\lambda_n+L,\lambda_2,\ldots,\lambda_{n-1},\lambda_1-L)
                       & for $i=0$ \\
               (\lambda_1,\ldots,\lambda_{i+1},\lambda_i,\ldots,\lambda_n)
                       & for $i=1,\ldots,n-1$\\}
\end{equation}
and $w\circ\lambda=w(\lambda+\bar{\rho})-\bar{\rho}$,
where $\bar\rho=(n-1,n-2,\ldots,1,0)$.

Let $\lambda\in \Pi_N^{(n,l)}$ and $\Lambda=\wt(\lambda)$.
If $\mu=w\circ\lambda$ is a partition, it is the unique one
such that $\wt(\mu)= w\cdot\Lambda$ and $|\mu|=|\lambda|$. So
the $\circ$ action corresponds to the projection of the
dot action on the classical part of the weight lattice.
If one forgets about columns of height $n$, one obtains
the true fusion algebra. In this context, the above algorithm
is due to Kac \cite{Kac} and Walton \cite{Wal1}.

The intersection of the orbit of a partition $\lambda\in\Pi^{(n,l)}$ 
under  $\W$ with the dominant
chamber admits a convenient graphical description.
It can be generated by sliding in the diagram of  $\lambda$
the  $L$-ribbons 
whose root lies in the first row (if there is no such ribbon,
then $\bs_\lambda\equiv 0$ unless $\lambda$ is a $(n,l)$-regular
$L$-core). 
The sign of an element of the orbit is the product of
the signs $(-1)^{h(R)-h(R')}$ for all ribbons $R$ of $\lambda$
where $R'$ is the new ribbon obtained by sliding $R$ 
and $h(R)$ denotes the heigth of $R$ \cite{GW,Cum}.

Let $r(R)$ be the position of the root of $R$, as in the picture
below. If $\mu$ is obtained from $\lambda$ by transforming $R$ into $R'$,
we define the $t$-equivalence by
\begin{equation}
s_\lambda \equiv_t (-1)^{h(R)-h(R')} t^{r(R')-r(R)} s_\mu \ .
\end{equation}
For example, with $n=3$ and $l=2$, $s_{644}$ is $t$-equivalent
to
$$
-t^{1} s_{743}, \ t^{2} s_{833}, \ t^{4} s_{77}, -t^{8} s_{11, 3}
, \ t^{10} s_{12, 2}
$$
as illustrated below:

\begin{center}
\footnotesize
\setlength{\unitlength}{0.00035in}
\begingroup\makeatletter\ifx\SetFigFont\undefined
\def\x#1#2#3#4#5#6#7\relax{\def\x{#1#2#3#4#5#6}}%
\expandafter\x\fmtname xxxxxx\relax \def\y{splain}%
\ifx\x\y   
\gdef\SetFigFont#1#2#3{%
  \ifnum #1<17\tiny\else \ifnum #1<20\small\else
  \ifnum #1<24\normalsize\else \ifnum #1<29\large\else
  \ifnum #1<34\Large\else \ifnum #1<41\LARGE\else
     \huge\fi\fi\fi\fi\fi\fi
  \csname #3\endcsname}%
\else
\gdef\SetFigFont#1#2#3{\begingroup
  \count@#1\relax \ifnum 25<\count@\count@25\fi
  \def\x{\endgroup\@setsize\SetFigFont{#2pt}}%
  \expandafter\x
    \csname \romannumeral\the\count@ pt\expandafter\endcsname
    \csname @\romannumeral\the\count@ pt\endcsname
  \csname #3\endcsname}%
\fi
\fi\endgroup
{\renewcommand{\dashlinestretch}{50}
\begin{picture}(15344,5239)(0,-10)
\path(22,5202)(1372,5202)(1372,3852)(922,3852)
\path(1372,5202)(1822,5202)(1822,4302)
	(2722,4302)(2722,3852)(1597,3852)
\path(1597,3852)(1372,3852)
\path(4072,5202)(4072,3852)(4972,3852)
	(4972,4752)(4072,4752)
\path(4072,5202)(5422,5202)(5422,3852)(4972,3852)
\path(5422,4752)(5872,4752)(5872,4302)
	(7222,4302)(7222,3852)(5422,3852)
\path(8572,5202)(8572,3852)(9472,3852)
	(9472,4752)(8572,4752)
\path(8572,5202)(9922,5202)(9922,3852)(9472,3852)
\path(9922,4302)(12172,4302)(12172,3852)(9922,3852)
\path(22,1152)(22,252)(922,252)
	(922,1152)(22,1152)
\path(922,1152)(1372,1152)(1372,702)
	(2722,702)(2722,252)(922,252)
\path(1372,1152)(3172,1152)(3172,252)
\path(3172,252)(2722,252)
\path(4072,1152)(4072,252)(4972,252)
	(4972,1152)(4072,1152)
\path(4972,1152)(5422,1152)(5422,702)
	(6772,702)(6772,252)(4972,252)
\path(6772,702)(9022,702)(9022,252)(6772,252)
\path(9922,1152)(9922,252)(10822,252)
	(10822,1152)(9922,1152)
\path(10822,702)(13072,702)(13072,252)(10822,252)
\path(13072,702)(15322,702)(15322,252)(13072,252)
\dottedline{45}(472,5202)(472,3852)
\dottedline{45}(922,5202)(922,4752)
\dottedline{45}(922,4752)(1822,4752)
\dottedline{45}(922,4302)(1822,4302)
\dottedline{45}(22,4302)(922,4302)
\dottedline{45}(1822,4302)(1822,3852)
\dottedline{45}(2272,4302)(2272,3852)
\dottedline{45}(4522,5202)(4522,3852)
\dottedline{45}(4972,5202)(4972,4752)
\dottedline{45}(5872,4302)(5872,3852)
\dottedline{45}(6322,4302)(6322,3852)
\dottedline{45}(6772,4302)(6772,3852)
\dottedline{45}(4972,4752)(5422,4752)
\dottedline{45}(4072,4302)(5872,4302)
\dottedline{45}(9022,5202)(9022,3852)
\dottedline{45}(9472,5202)(9472,4752)
\dottedline{45}(10372,4302)(10372,3852)
\dottedline{45}(10822,4302)(10822,3852)
\dottedline{45}(11272,4302)(11272,3852)
\dottedline{45}(11722,4302)(11722,3852)
\dottedline{45}(9472,4752)(9922,4752)
\dottedline{45}(8572,4302)(9922,4302)
\dottedline{45}(472,1152)(472,252)
\dottedline{45}(1372,702)(1372,252)
\dottedline{45}(1822,1152)(1822,252)
\dottedline{45}(2272,1152)(2272,252)
\dottedline{45}(2722,1152)(2722,702)
\path(22,5202)(22,3852)(922,3852)
	(922,4752)(22,4752)
\dottedline{45}(22,702)(1372,702)
\put(14872,27){\makebox(0,0)[lb]{\smash{{{\SetFigFont{7}{14.4}{rm} 12}}}}}
\dottedline{45}(2722,702)(3172,702)
\dottedline{45}(4522,1152)(4522,252)
\dottedline{45}(4072,702)(4972,702)
\dottedline{45}(4972,702)(5422,702)
\dottedline{45}(5422,702)(5422,252)
\dottedline{45}(5872,702)(5872,252)
\dottedline{45}(6322,702)(6322,252)
\dottedline{45}(7222,702)(7222,252)
\dottedline{45}(7672,702)(7672,252)
\dottedline{45}(8122,702)(8122,252)
\dottedline{45}(8572,702)(8572,252)
\dottedline{45}(9922,702)(10822,702)
\dottedline{45}(10372,1152)(10372,252)
\dottedline{45}(11272,702)(11272,252)
\dottedline{45}(11722,702)(11722,252)
\dottedline{45}(12172,702)(12172,252)
\dottedline{45}(12622,702)(12622,252)
\dottedline{45}(13522,702)(13522,252)
\dottedline{45}(13972,702)(13972,252)
\dottedline{45}(14422,702)(14422,252)
\dottedline{45}(14872,702)(14872,252)
\path(1597,4977)(1597,4077)(2497,4077)
\path(247,4977)(1147,4977)(1147,4077)
\path(4297,4977)(5197,4977)(5197,4077)
\path(5647,4527)(5647,4077)(6997,4077)
\path(8797,4977)(9697,4977)(9697,4077)
\path(10147,4077)(11947,4077)
\path(1147,927)(1147,477)(2497,477)
\path(1597,927)(2947,927)(2947,477)
\path(5197,927)(5197,477)(6547,477)
\path(6997,477)(8797,477)
\path(11047,477)(12847,477)
\path(13297,477)(15097,477)
\put(247,3627){\makebox(0,0)[lb]{\smash{{{\SetFigFont{7}{14.4}{rm}1}}}}}
\put(697,3627){\makebox(0,0)[lb]{\smash{{{\SetFigFont{7}{14.4}{rm}2}}}}}
\put(1147,3627){\makebox(0,0)[lb]{\smash{{{\SetFigFont{7}{14.4}{rm}3}}}}}
\put(1597,3627){\makebox(0,0)[lb]{\smash{{{\SetFigFont{7}{14.4}{rm}4}}}}}
\put(2047,3627){\makebox(0,0)[lb]{\smash{{{\SetFigFont{7}{14.4}{rm}5}}}}}
\put(2497,3627){\makebox(0,0)[lb]{\smash{{{\SetFigFont{7}{14.4}{rm}6}}}}}
\put(6997,3627){\makebox(0,0)[lb]{\smash{{{\SetFigFont{7}{14.4}{rm}7}}}}}
\put(11947,3627){\makebox(0,0)[lb]{\smash{{{\SetFigFont{7}{14.4}{rm}8}}}}}
\put(2497,27){\makebox(0,0)[lb]{\smash{{{\SetFigFont{7}{14.4}{rm}6}}}}}
\put(2947,27){\makebox(0,0)[lb]{\smash{{{\SetFigFont{7}{14.4}{rm}7}}}}}
\put(6547,27){\makebox(0,0)[lb]{\smash{{{\SetFigFont{7}{14.4}{rm}6}}}}}
\put(8572,27){\makebox(0,0)[lb]{\smash{{{\SetFigFont{7}{14.4}{rm}   11}}}}}
\put(12847,27){\makebox(0,0)[lb]{\smash{{{\SetFigFont{7}{14.4}{rm}9}}}}}
\end{picture}
}

\end{center}
\normalsize

The power of $t$ is just the opposite of the
coefficient of the null root $\delta$
in $w\cdot\Lambda$, where $\Lambda=\wt(\lambda)$.
In our example, we start from $\lambda=(6,4,4)$. We have
$n=3$ and $L=5$ so that we use weights of level $5-3=2$.
Therefore, $\wt(\Lambda)=2\Lambda_1$. The orbit of $\Lambda$ is
\begin{eqnarray*}
\Lambda=2\Lambda_1 & \leftrightarrow & s_{644} \\
s_0\cdot \Lambda = -2\Lambda_0+3\Lambda_1+\Lambda_2-\delta
                & \leftrightarrow & -t s_{743}\\
s_0s_2\cdot\Lambda = -3\Lambda_0+5\Lambda_1-2\delta
                  & \leftrightarrow & t^2 s_{833} \\
s_0s_1\cdot\Lambda = -5\Lambda_0+7\Lambda_2-4\delta 
                & \leftrightarrow & t^4s_{77} \\
s_0s_2s_1\cdot\Lambda = -8\Lambda_0+7\Lambda_1+3\Lambda_2-8\delta
                 & \leftrightarrow & -t^8s_{11,3} \\
s_0s_2s_1s_0\cdot\Lambda = -10\Lambda_0+10\Lambda_1+2\Lambda_2-10\delta
                   & \leftrightarrow & t^{10} s_{12,2} \ .
\end{eqnarray*}

The $q$-fusion coefficients are defined by applying this
$t$-reduction algorithm to the $q$-analogues of unrestricted
products of Schur functions, with $t=q$ for the spin/charge/energy
$q$-analogues, and
$t=1/q$ for the cospin/cocharge $q$-analogues.

For example, the cospin $q$-analogue of the cube $(s_{21})^3$
of the Schur function $s_{21}$ (restricted to partitions of length
$\le 4$) is equal to
\begin{eqnarray*}\fl
\tilde G_{666333}^{(3)} = &
s_{63}+(q+q^2)s_{621}+q^3s_{6111} +(q+q^2)s_{54} \\
&+ (q+2q^2+2q^3+q^4)s_{531} + (q^2+2q^3+q^4)s_{522}
+(q^2+2q^3+2q^4+q^5)s_{5211}\\
&+ (q^2+2q^3+q^4)s_{441} + (q^2+2q^3+3q^4+2q^5)s_{432}
+(2q^3+3q^4+3q^5+q^6)s_{4311}\\
&+(q^3+3q^4+3q^5+2q^6)s_{4221} +(q^3+q^6)s_{333} +
(2q^4+3q^5+2q^6+q^7)s_{3321}\\
&+ (q^5+2q^6+q^7)s_{3222} \,.
\end{eqnarray*}
Now, take $n=4$ and $l=2$, so that only the partitions
$(3321)$ and $(3222)$ are $(4,2)$-restricted. The $q$-reduction
algorithm produces as $q$-analogue of $(\bar s_{21})^3$ in
${\cal F}^{(4,2)}$
\begin{equation}
q^7 \bar s_{3221} + (q^6+q^7) \bar s_{3222} \,.
\end{equation}

We conjecture that our restricted $q$-analogues are always positive,
and that
for a product of rectangular Schur functions, the spin $q$-analogues
coincide (up to an overall power of $q$) with the ones defined
by equation (\ref{LRpath}). We can prove this 
from the results of \cite{KMOTU} in the case where all the
$\mu^{(i)}$ are equal to the same row (or column) partition.

The modified Hall-Littlewood function $Q'_\mu$ is a
$q$-analogue of the product $h_{\mu_1}h_{\mu_2}\ldots h_{\mu_r}$.
In this case, our conjecture states that the restricted Kostka
polynomial $\bar K_{\lambda\mu}(q)$ (coefficient of $\bar s_\lambda$
in $\bar Q'_\mu$) is equal to the
restricted 1d-configuration sum $X^{(l-k)}_\mu(\lambda)$ of \cite{HKKOTY}.
For example, in Example 3.1 of \cite{HKKOTY}  it is found that
$X_{2211}^{(2)}(321)=q$.
On the other hand, the above algorithm gives for the restricted
Kostka polynomial
\begin{eqnarray*}\fl
\bar K_{321,2211}(q)&=&
K_{321,2211}(q)-q^{-1}K_{42,2211}(q)+q^{-3}K_{6,2211}(q)
\\
&=& (q+2q^2+q^3) -q^{-1}(2q^3+q^4+q^5)+q^{-3}(q^7) = q\,.
\end{eqnarray*}

\section{Examples}

\subsection{Virasoro characters}

We consider the $\slchap_2$ case ($L=l+2$). Assume $\mu^{(1)}=\cdots
=\mu^{(r)}=(1)$. We set
\[
\bk_{(N+b,N)}(q) = 
\bar c^{(N+b,N)}_{(1),\ldots,(1)}(q)=\bar K_{(N+b,N),(1^{2N+b})}(q).
\]
The intersection of the orbit of a dominant integral weight of $A^{(1)}_1$
under the affine Weyl group with the dominant
chamber can be explicitly computed.
As a result, we obtain an explicit formula for
the restricted Kostka polynomial 
$$
\bk_{(N+b,N)}(q) =
\sum_{n\ge 0} q^{-n(Ln+b+1)}k_{(N+Ln+b,N-Ln)}(q) 
- \sum_{n\ge 1} q^{-n(Ln-b-1)}k_{(N+Ln-1,N-Ln+b+1)}(q) \ .
$$
Taking into account the expression
$$
k_{(\lambda_1,\lambda_2)}(q)=
q^{{\lambda_1\choose 2} + {\lambda_2\choose 2}}
\qbin{\lambda_1+\lambda_2}{\lambda_1}
-
q^{{\lambda_1+1 \choose 2} +{\lambda_2-1\choose 2}}
\qbin{\lambda_1 +\lambda_2}{\lambda_1+1} \ ,
$$
one obtains
\begin{eqnarray}\fl
\bk_{(N+b,N)}(q) =
q^{N(N+b-1) +b(b-1)/2}
\left(
\sum_{n=-\infty}^\infty
q^{n[L(L-1)n+Lb-b-1]}\qbin{2N+b}{N-Ln} \right. \nonumber \\ 
\left. -
\sum_{n=-\infty}^\infty
q^{n[L(L-1)n+Lb+2L-b-1]   +b+1    }\qbin{2N+b}{N-Ln-1}
\right) \,.
\end{eqnarray}

A restricted path in $\P^{(N+b,N)}_{(1),\ldots,(1)}$ can be interpreted as 
an oriented path on $P^+_l$ starting from $l\Lambda_0$ and arriving at
$(l-b)\Lambda_0+b\Lambda_1$.
More generally, paths of length $2N+b-a$ starting from 
$(l-a)\Lambda_0+a\Lambda_1$ and arriving to $(l-b)\Lambda_0+b\Lambda_1$ 
are encoded by the skew diagram
of shape $(N+b,N)/(a)$. The reduced skew Kostka polynomials
are obtained in a similar way.
\begin{eqnarray}\fl
\bk_{(N+b,N)/(a)}(q) =
q^{N(N+b-a-1)+(b-a)(b-a-1)/2}
\sum_{n=-\infty}^\infty
\left(
q^{n[L(L-1)n+(L-1)b-La-1]} \qbin{2N+b-a}{N-Ln} \right. \nonumber\\
-  \left. \sum_{n=-\infty}^\infty
q^{[(L-1)n+a+1][Ln+b+1]}\qbin{2N+b-a}{N+Ln+b+1} \right) \ .
\end{eqnarray}
In the limit $N\rightarrow \infty$, we get 
\begin{eqnarray}
\lim_{N\rightarrow\infty} q^{-N(N+b-a-1)-(b-a)(b-a-1)/2}
\bk_{(N+b,N)/(a)}(q)
\nonumber \\
= {1\over\varphi(q)} \sum_{n=-\infty}^\infty
\left( q^{n[L(L-1)n+(L-1)b-La-1]}
-q^{[(L-1)n+a+1][Ln+b+1]} \right).
\end{eqnarray}

We shall relate this series to a Virasoro charcter.
We use the following conventions for the Virasoro characters.
The unitary series  minimal models are parametrized by 
three integers giving their
central charges
\begin{equation}
c = 1 - {6\over m(m+1)} \qquad (m\in\N,\ m\ge 2)
\end{equation}
and conformal weights
\begin{equation}
h=h_{r,s} = {[(m+1)r-ms]^2-1\over 4m(m+1)}\, \quad r=1,\ldots,m-1, \
s=1,\ldots,r \ .
\end{equation}
The character of 
the module  $M^{(m)}(r,s)=L(c,h)$ with $c$ and $h$ as above
is given by the Rocha-Caridi formula \cite{RC}
\begin{equation}
\chi^{(m)}_{r,s}(q) = {1\over\varphi(q)} \sum_{n=-\infty}^\infty
\left(q^{a_-(n)}-q^{a_+(n)} \right)
\end{equation}
where
\begin{equation}
a_\pm(n)=
{[2m(m+1)n+(m+1)r\pm ms]^2-1 \over 4m(m+1)} \,.
\end{equation}
Taking $m=L-1$, $r=a+1$, $s=b+1$,
we find
\begin{equation}
\lim_{N\rightarrow\infty} q^{-N(N+b-a-1)-(b-a)(b-a-1)/2}
\bk_{(N+b,N)/(a)}(q)
=
q^{-h} \chi^{(L-1)}_{a+1,b+1}(q) \, .
\end{equation}

\subsection{Nonabelian theta functions}

Let ${\cal SU}_\Sigma(n)$ be the moduli space of semi-stable holomorphic
vector bundles with trivial
determinant over a compact Riemann surface $\Sigma$ of genus $g$ 
({\it cf. } \cite{BL}). All line bundles on ${\cal SU}_\Sigma(n)$
are powers of the so-called determinant bundle ${\cal L}$.
The sections of ${\cal L}^l$ are called nonabelian theta functions of
rank $n$ and level $l$ on $\Sigma$. The dimension $h^0({\cal L}^l)$
of $H^0({\cal SU}_\Sigma(n),{\cal L}^l)$ is obtained by a calculation
in the fusion algebra ${\cal F}^{(n,l)}$ as follows \cite{Ver,BL}.
For a partition $\lambda$ interpreted as a dominant weight of $sl_n$,
let $\lambda^*=-w_0(\lambda)$ be the highest weight of the dual
representation $(V_\lambda)^*$. Now, let 
\begin{equation}
\omega=\sum_{\lambda\subset(l^{n-1})} \bar s _\lambda \bar s_{\lambda^*}
\ \ \in \ \ {\cal F}^{(n,l)} \,.
\end{equation}
Then, $h^0({\cal L}^l)$ is equal to the constant term in $\omega^g$.
For example, with $n=2$, $h^0({\cal L})=2^g$,
$h^0({\cal L}^2)=2^{g-1}(2^g+1)$, and so on.
If we compute $\omega^g$ by using the cospin $q$-analogues of the
fusion coefficients, we obtain $q$-analogues $h_q^0({\cal L}^l)$
of these numbers.
With $n=2$ as above, one finds $h_q^0({\cal L})=(1+q)^g$
and $h_q^0({\cal L}^2)={1\over 2}[(1+q^2)^g+(1+q)^{2g}]$.
This suggests the existence of some natural filtration
of the space of nonabelian theta functions.

\end{document}